\documentclass[12pt]{article}
\usepackage[centertags]{amsmath}
\usepackage{amsfonts}
\usepackage{amssymb}
\usepackage{amsthm}
\usepackage{amsmath,mathrsfs}
\usepackage{amsmath,amscd}
\title{\textbf{A Stringy Generalization of the Kontsevich Integral}}
\author{Renaud Gauthier \footnote{rg.mathematics@gmail.com}\\ \\Lycee Albert Camus}
\theoremstyle{definition}
\newtheorem{ceedee}{Definition}[section]

\newcommand{\beq}{\begin{equation}}
\newcommand{\eeq}{\end{equation}}

\newcommand{\dlog}{\text{dlog}}

\newcommand{\barA}{\overline{\mathcal{A}}}
\begin{document}
\maketitle
\begin{abstract}
We introduce a ``minimal" Kontsevich integral that generates the original Kontsevich integral while at the same time producing ribbons whose boundaries are the braids on which the minimal Kontsevich integral is evaluated. We generalize the definition of the Kontsevich integral to that of graphs in $\mathbb{C} \times I$ and study the behavior of such expressions as different graphs are brought together, thus leading to a 2-dimensional generalization of the Kontsevich integral.
\end{abstract}
\newpage

\section{The Kontsevich Integral}
The Kontsevich integral \cite{K} is a functional on knots that can be seen as a generalization of the Gauss integral. It is a graded sum of chord diagrams times some coefficients that are essentially integrals of powers of log differentials. Chord diagrams are the knots used as an argument in $Z$ with horizontal dashed chords stretching between its strands. The Kontsevich integral $Z$ is highly dependent on a choice of time axis. Further, if it is invariant under horizontal deformations that keep the local extrema of knots used as an argument fixed, it is not invariant under translations for which such extrema are moved. Thus $Z$ depends on a path along which graphs are translated, as well as rotations, thus presenting the Kontsevich integral as a map of what appears to be Riemann surfaces, but the exact structure of which we will study in this paper.\\

Before introducing this integral, we define the algebra $\mathcal{A}$ ~\cite{K} in which it takes its values. For a singular oriented knot whose only singularities are transversal self-intersections, the preimage of each singular crossing under the embedding map defining the knot yields a pair of distinct points on $S^1$. Each singular point in the image therefore yields a pair of points on $S^1$ that are conventionally connected by a chord for book keeping purposes. A knot with $m$ singular points will yield $m$ distinct chords on $S^1$. One refers to such a circle with $m$ chords on it as a chord diagram of degree $m$, the degree being the number of chords. The support of the graph is an oriented $S^1$, and it is regarded up to orientation preserving diffeomorphisms of the circle. More generally, for a singular oriented link all of whose singularities are double-crossings, preimages of each singular crossing under the embedding map defining the link yield pairs of distinct points on possibly different circles depending on whether the double crossing was on a same component or between different components of the link. One also connects points making a pair by a chord. A link with $m$ singular points will yield $m$ chords on $\coprod S^1$. One calls such a graph a chord diagram. The support is $\coprod S^1$ regarded up to orientation preserving diffeomorphism of each $S^1$.\\

One denotes by $\mathcal{D}$ the complex vector space spanned by chord diagrams with support $S^1$. There is a grading on $\mathcal{D}$ given by the number of chords featured in a diagram. If $\mathcal{D}^{(m)}$ denotes the subspace of chord diagrams of degree $m$, then one writes:
\beq
\mathcal{D}=\oplus_{m\geq 0} \mathcal{D}^{(m)}
\eeq
One quotients this space by the 4-T relation which locally looks like:\\
\setlength{\unitlength}{1cm}
\begin{picture}(1,2)(-1,-0.5)
\multiput(0,0.75)(0.2,0){5}{\line(1,0){0.1}}
\multiput(1,0.25)(0.2,0){5}{\line(1,0){0.1}}
\put(0,0.9){\vector(0,1){0.2}}
\put(1,0.9){\vector(0,1){0.2}}
\put(2,0.9){\vector(0,1){0.2}}
\linethickness{0.3mm}
\put(0,0){\line(0,1){1}}
\put(1,0){\line(0,1){1}}
\put(2,0){\line(0,1){1}}
\put(2.5,0.5){$+$}
\end{picture}
\setlength{\unitlength}{1cm}
\begin{picture}(1,2)(-3,-0.5)
\multiput(0,0.75)(0.2,0){5}{\line(1,0){0.1}}
\multiput(0,0.25)(0.2,0){10}{\line(1,0){0.1}}
\put(0,0.9){\vector(0,1){0.2}}
\put(1,0.9){\vector(0,1){0.2}}
\put(2,0.9){\vector(0,1){0.2}}
\linethickness{0.3mm}
\put(0,0){\line(0,1){1}}
\put(1,0){\line(0,1){1}}
\put(2,0){\line(0,1){1}}
\put(2.5,0.5){$=$}
\end{picture}
\setlength{\unitlength}{1cm}
\begin{picture}(1,2)(-5,-0.5)
\multiput(0,0.25)(0.2,0){5}{\line(1,0){0.1}}
\multiput(1,0.75)(0.2,0){5}{\line(1,0){0.1}}
\put(0,0.9){\vector(0,1){0.2}}
\put(1,0.9){\vector(0,1){0.2}}
\put(2,0.9){\vector(0,1){0.2}}
\linethickness{0.3mm}
\put(0,0){\line(0,1){1}}
\put(1,0){\line(0,1){1}}
\put(2,0){\line(0,1){1}}
\put(2.5,0.5){$+$}
\end{picture}
\setlength{\unitlength}{1cm}
\begin{picture}(1,2)(-7,-0.5)
\multiput(0,0.25)(0.2,0){5}{\line(1,0){0.1}}
\multiput(0,0.75)(0.2,0){10}{\line(1,0){0.1}}
\put(0,0.9){\vector(0,1){0.2}}
\put(1,0.9){\vector(0,1){0.2}}
\put(2,0.9){\vector(0,1){0.2}}
\linethickness{0.3mm}
\put(0,0){\line(0,1){1}}
\put(1,0){\line(0,1){1}}
\put(2,0){\line(0,1){1}}
\end{picture}\\ \\
where solid lines are intervals on $S^1$ on which a chord foot rests, and arrows indicate the orientation of each strand. One further quotients this space by the framing independence relation: if a chord diagram has a chord forming an arc on $S^1$ with no other chord ending in between its feet, then the chord diagram is set to zero. The resulting quotient space is the complex vector space generated by chord diagrams mod the 4-T relation and framing independence and is denoted by $\mathcal{A}$. The grading of $\mathcal{D}$ is preserved by the quotient, inducing a grading on $\mathcal{A}$:
\beq
\mathcal{A}=\oplus_{m \geq 0}\mathcal{A}^{(m)}
\eeq
where $\mathcal{A}^{(m)}$ is obtained from $\mathcal{D}^{(m)}$ upon modding out by the 4-T and the framing independence relations. All this carries over to the case of links by formally extending the 4-T relation to the case of $q$ disjoint copies of the circle in the case of a $q$-components link, and the resulting $\mathbb{C}$-vector space will be denoted $\mathcal{A}(\coprod_q S^1)$.\\

The connected sum of circles can be extended to chorded circles, thereby defining a product on $\mathcal{A}$, making it into an associative and commutative algebra. The Kontsevich integral will be valued in the graded completion $\overline{\mathcal{A}}=\prod_{m \geq 0}\mathcal{A}^{(m)}$ of the algebra $\mathcal{A}$.\\

As far as knots are concerned, one works with Morse knots, geometric tangles and graphs whose vertices are curved lines. We distinguish graphs that are initially given in the argument of $Z$ from those that result from the gluing of two distinct graphs. The reason for this distinction is that initial graphs will be univalent, trivalent (y or $\lambda$-shaped) or 4-valent (X-shaped) so that their corresponding Kontsevich integral is non-singular. However graphs that result from the gluing of two graphs may have vertices that do not fall in either of these categories and thus may likely result in the Kontsevich integral of such graph being singular as we will see later. Having said that, one considers all such geometric pictures being embedded in $\mathbb{R}^3$ a decomposition of which can be given as the product of the complex plane and the real line: $\mathbb{R}^3=\mathbb{R}^2\times \mathbb{R} \simeq \mathbb{C}\times \mathbb{R}$, with local coordinates $z$ on the complex plane and $t$ on the real line for time. A morse knot $K$ is such that $t\circ K$ is a Morse function on $S^1$. An acceptable graph for our purposes is an initial graph as defined above, so that after a possible rotation none of its edges end up being a horizontal edge, something that could happen should one of its edges be a straight line. If one denotes by $Z$ the Kontsevich integral functional on knots, if $K$ is a Morse knot, one defines ~\cite{K}, ~\cite{ChDu}:
\beq
Z(K):=\sum_{m\geq 0} \frac{1}{(2 \pi i)^m}\int_{t_{min}< t_1<...<t_m<t_{max}}\sum_{P\; applicable}(-1)^{\varepsilon(P)}D_P\prod_{1\leq i \leq m}\dlog \vartriangle \!\!z[P_i] \label{IK}
\eeq
where $t_{min}$ and $t_{max}$ are the min and max values of $t$ on $K$ respectively, $P$ is an $m$-vector each entry of which is a pair of points on the image of the knot $K$, $P=(P_1,...,P_m)$, where the $i$-th entry $P_i$ corresponds to a pair of points on the knot. One refers to such $P$'s as pairings. If one further situates these paired points at some height $t_i$, and denote these two points by $z_i$ and $z'_i$, then we define $\vartriangle \!\! z[P_i]:=z_i-z'_i$. One denotes by $K_P$ the knot $K$ with $m$ pairs of points placed on it following the prescription given by $P$, along with chords connecting such points at a same height. A pairing is said to be applicable if each entry is a pair of two distinct points on the knot, at the same height ~\cite{ChDu}. We will assume that all chords are horizontal on knots and will drop the adjective applicable, simply referring to $P$'s as pairings. One denotes by $\varepsilon(P)$ the number of those points ending on portions of $K$ that are locally oriented down. For example if $P=(z(t),z'(t))$ and $K$ is decreasing at $z(t)$, then it will contribute 1 to $\varepsilon(P)$. One also define the length of $P$ to be $|P|$, the number of pairings it is a combination of. If we denote by $\iota_K$ the embedding defining the knot $K$ then $D_P$ is defined to be the chord diagram one obtains by taking the inverse image of $K_P$ under $\iota_K$, that is $D_P=\iota_K^{-1} K_P$. This generalizes immediately to the case of Morse links, and in this case the geometric coefficient will not be an element of $\barA$ but will be an element of $\barA (\coprod_q S^1)$ if the argument of $Z$ is a $q$-components link. Observe that $Z(L) \in \overline{\mathcal{A}}(\amalg_q S^1)$ is known once $L_P$ is known for all $P$'s. This is what is referred to as a tangle chord diagram \cite{ChDu} and sometimes the Kontsevich integral is given not as an element of $\overline{\mathcal{A}}(\amalg S^1)$ but as an element of $\overline{\mathcal{A}}(L)$ and is written instead exactly as in \eqref{IK} except that instead of using chord diagrams $D_P$ tangle chord diagrams $L_P$ are used. This generalizes to the case of a geometric braid or even a graph $\Gamma$ by using $\Gamma_P$ instead of $L_P$.\\

In section 2, we introduce the configuration space of $N$ unordered points in the complex plane. In section 3 we present the general notion of chord diagrams. In section 4 we present a minimal Kontsevich integral that we regard as an equivalence between objects in $\overline{\mathcal{A}}$ and links, but which can also be used to generate the original Kontsevich integral. In section 5 we show that we can make the Kontsevich integral more dynamic by making it time dependent as well as dependent on rotations, thus introducing cylinders in the picture on which $Z$ is defined.\\

\section{The configuration space of $N$ points in the plane}
A link in $S^3$ is ambient isotopic to a closed braid ~\cite{A} ~\cite{JB}, so that one can deform a link into a braid part, outside of which all its strands are parallel. For a given link, let $N$ be the number of strands of its braid part. $N$ will depend on the link we have chosen. The transversal intersection of these $N$ strands with the complex plane will yield a set of $N$ distinct points, each point resulting from the intersection of one strand with this plane. It is natural then to study, for any given $N$, the space $X_N$ defined as the configuration space of $N$ distinct unordered points in the complex plane:
\beq
X_N:=\{(z_1,...,z_n) \in \mathbb{C}^N | z_i=z_j \Rightarrow i=j\}/ S_N=(\mathbb{C}^N-\Delta)/S_N
\eeq
where $S_N$ is the permutation group on $N$ elements and $\Delta$ is the big diagonal in $\mathbb{C}^N$. The labeling of points of $X_N$ is not induced by any ordering on $\mathbb{C}^N$ but rather is a way to locate the $N$ points in the complex plane whose collection defines a single point of $X_N$. We will sometimes write $\sum_{1 \leq i \leq N}[z_i]$ instead of $\{z_1,...,z_N\}$ to represent points in configuration space. The points $z_1,...,z_N$ of the complex plane defining a point $Z=\sum_{1 \leq i \leq N}[z_i]$ of $X_N$ will be referred to as the $N$ defining points of $Z$. We consider the topology $\tau$ on $X_N$ generated by open sets of the form $U=\{U_1,...,U_N\}$ where the $U_i,\,1 \leq i \leq N$ are non-overlapping open sets in the complex plane. We will also refer to those open sets $U_1,...,U_N$ as the $N$ defining open sets of the open set $U$ of $X_N$. \\

We review the basic terminology pertaining to braids as presented in ~\cite{JB} since we will work with braids in what follows. The pure braid group of $\mathbb{C}^N$ is defined to be $\pi_1 (\mathbb{C}^N - \Delta)$, and the braid group of $\mathbb{C}^N$ is defined to be $\pi_1 (X_N)$. A braid is an element of this latter group. If $q$ denotes the regular projection map from $\mathbb{C}^N - \Delta$ to $X_N$, $Z=(z_1,...,z_N) \in \mathbb{C}^N - \Delta$, $qZ \in X_N$, then $\gamma \in \pi_1 (X_N,qZ)$ based at $qZ$ is given by a loop $\gamma =\{\gamma_1,...,\gamma_N \}$ which lifts uniquely to a path in $\mathbb{C}^N - \Delta$ based at $Z$ that without loss of generality we will denote by the same letter $\gamma$. Then we have $\gamma =(\gamma_1 ,...,\gamma_N )$. The graph of the $i$-th coordinate of $\gamma$ is defined to be $\Gamma_i := \{(\gamma_i (t),t)\;|\;t \in I \}$, $1 \leq i \leq N$. Each such graph $\Gamma_i$ defines an arc $\tilde{\gamma}_i \in \mathbb{C} \times I$ and $\tilde{\gamma}:=\cup_{1 \leq i \leq N} \tilde{\gamma}_i \in \mathbb{C}\times I $ is called a geometric braid, which we will refer to as the lift of $\gamma$. As such it is open, and its closure is a closed braid.

\section{Chord diagrams}
We will be interested in considering chord diagrams with support graphs in $\mathbb{C} \times I$, so for that purpose one considers a more general definition of chord diagrams than the one presented in the introduction which was sufficient to discuss the Kontsevich integral of knots.
\begin{ceedee}[\cite{LM}]
Let $X$ be a one dimensional, compact, oriented, smooth manifold with corners with numbered components. A chord diagram with support on $X$ is a set of finitely many unordered pairs of distinct non-boundary points on $X$ defined modulo orientation and component preserving homeomorphisms. One realizes each pair geometrically by drawing a dashed line, or chord, stretching from one point to the other. One denotes by $\mathcal{A}(X)$ the $\mathbb{C}$-vector space spanned by chord diagrams with support on $X$ modulo the framing indepence relation as well as the 4-T relation: if $i$, $j$ and $k$ are indices for components of $X$ on which chords are ending, then locally the 4-T relation can be written:

\setlength{\unitlength}{1cm}
\begin{picture}(1,2)(0,-0.5)
\multiput(0,0.75)(0.2,0){5}{\line(1,0){0.1}}
\multiput(1,0.25)(0.2,0){5}{\line(1,0){0.1}}
\put(0,0.9){\vector(0,1){0.2}}
\put(1,0.9){\vector(0,1){0.2}}
\put(2,0.9){\vector(0,1){0.2}}
\linethickness{0.3mm}
\put(0,0){\line(0,1){1}}
\put(1,0){\line(0,1){1}}
\put(2,0){\line(0,1){1}}
\put(0.1,-0.2){$i$}
\put(1.1,-0.2){$j$}
\put(2.1,-0.2){$k$}
\put(2.5,0.5){$+$}
\end{picture}
\setlength{\unitlength}{1cm}
\begin{picture}(1,2)(-2,-0.5)
\multiput(0,0.75)(0.2,0){5}{\line(1,0){0.1}}
\multiput(0,0.25)(0.2,0){10}{\line(1,0){0.1}}
\put(0,0.9){\vector(0,1){0.2}}
\put(1,0.9){\vector(0,1){0.2}}
\put(2,0.9){\vector(0,1){0.2}}
\linethickness{0.3mm}
\put(0,0){\line(0,1){1}}
\put(1,0){\line(0,1){1}}
\put(2,0){\line(0,1){1}}
\put(0.1,-0.2){$i$}
\put(1.1,-0.2){$j$}
\put(2.1,-0.2){$k$}
\put(2.5,0.5){$=$}
\end{picture}
\setlength{\unitlength}{1cm}
\begin{picture}(1,2)(-4,-0.5)
\multiput(0,0.25)(0.2,0){5}{\line(1,0){0.1}}
\multiput(1,0.75)(0.2,0){5}{\line(1,0){0.1}}
\put(0,0.9){\vector(0,1){0.2}}
\put(1,0.9){\vector(0,1){0.2}}
\put(2,0.9){\vector(0,1){0.2}}
\linethickness{0.3mm}
\put(0,0){\line(0,1){1}}
\put(1,0){\line(0,1){1}}
\put(2,0){\line(0,1){1}}
\put(0.1,-0.2){$i$}
\put(1.1,-0.2){$j$}
\put(2.1,-0.2){$k$}
\put(2.5,0.5){$+$}
\end{picture}
\setlength{\unitlength}{1cm}
\begin{picture}(1,2)(-6,-0.5)
\multiput(0,0.25)(0.2,0){5}{\line(1,0){0.1}}
\multiput(0,0.75)(0.2,0){10}{\line(1,0){0.1}}
\put(0,0.9){\vector(0,1){0.2}}
\put(1,0.9){\vector(0,1){0.2}}
\put(2,0.9){\vector(0,1){0.2}}
\linethickness{0.3mm}
\put(0,0){\line(0,1){1}}
\put(1,0){\line(0,1){1}}
\put(2,0){\line(0,1){1}}
\put(0.1,-0.2){$i$}
\put(1.1,-0.2){$j$}
\put(2.1,-0.2){$k$}
\end{picture}

One defines the degree of a chord diagram to be the number of chords a chord diagram has, and we call it the chord degree. This induces a graded decomposition of the space $\mathcal{A}(X)$:
\beq
\mathcal{A}(X)=\bigoplus_{m\geq 0}\mathcal{A}^{(m)}(X)
\eeq
where $\mathcal{A}^{(m)}(X)$ is the $\mathbb{C}$-vector space of chord diagrams of degree $m$ with support on $X$. One writes $\overline{\mathcal{A}}(X)$ for the graded completion of $\mathcal{A}(X)$.
\end{ceedee}

We will initially be interested in the case where $X$ is a geometric braid $\tilde{\gamma} \in \mathbb{C} \times I$ corresponding to some loop $\gamma$ in $X_N$. The strands are oriented up, $t=0$ being the bottom plane of the space $\mathbb{C}\times I$ in which the braid is embedded, $t=1$ corresponding to the top plane. Since indices for pairings match those for the times at which they are located, chords will be ordered from the bottom up. For $m=1$, a chord will stretch between two strands, say the strands indexed by $i$ and $j$, and we will denote such a chord diagram by  $|ij\rangle \in \mathcal{A}(\tilde{\gamma})$, corresponding to the  pairing $(ij)$ in this case. If we want to insist that the skeleton of the chord diagram is a given geometric braid $\tilde{\gamma}$ then we write $|ij\rangle (\tilde{\gamma})$. In certain situations it will be necessary to also indicate at which point along the braid is the chord situated for location purposes. Once we have $|ij\rangle (\tilde{\gamma})$, it is sufficient to have the height $t \in I$ at which we have to place the chord $|ij\rangle$ on $\tilde{\gamma}$ and $|ij\rangle (\tilde{\gamma})(t)$ is defined to be a chord between the $i$-th and $j$-th strands of $\tilde{\gamma}$ at height $t$, or equivalently a chord between $(\gamma_i (t),t)$ and $(\gamma_j (t),t)$. In that case we work with a representative of the class defining the chord diagram $|ij\rangle (\tilde{\gamma})$. \\

For the purpose of reconstructing links from chord diagrams, we will be interested in chord diagrams supported at a point of $X_N$. For a point $Z=\{z_1,...,z_N\} \in X_N$, some $P=(k,l)$, $1 \leq k \neq l \leq N$, $|P\rangle(Z) \in \mathcal{A}(Z)$ is a chord between $z_k$ and $z_l$ in $X_N$. We denote by $\mathcal{A}(X_N)$ the complex vector space spanned by all such elements, and by $\overline{\mathcal{A}}(X_N)$ its graded completion.\\

We will also be interested in working with elements of $\mathcal{A}^{(1)}(X) \otimes \Omega^1 (\log  \mathbb{C})$ with $X$ to be determined, that we denote by $|ij\rangle \dlog (z_i-z_j)$. In this notation if $\tilde{\gamma} \in \mathbb{C} \times I$ is a geometric braid obtained from lifting a loop $\gamma$ in $X_N$, if we arbitrarily index the $N$ strands of $\tilde{\gamma}$, then the $k$-th strand is obtained from lifting a path in the complex plane given by some function $z(t)$, $t \in I$. For a chord $|ij\rangle$ between the $i$-th and the $j$-th strands which are the respective lifts of paths $\gamma_i$ and $\gamma_j$ in the complex plane given by functions $z_i(t)$ and $z_j(t)$, $t\in I$, then $z_i-z_j$ is the difference of two such functions. This leads us to defining the subspace $\Omega^1 (\log \!\vartriangle \! \!\mathbb{C})$ of log differential functionals on $\mathbb{C}$, defined by $\dlog(\vartriangle \!\! z[z_1, z_2])=\dlog(z_1-z_2)$. We have a projection:
\begin{align}
\Omega^1(log \!\vartriangle \!\!\mathbb{C}) & \xrightarrow{p_2} (\mathbb{C}^2-\Delta)/S_2 \\
dlog(z_i-z_j) &\mapsto \{z_i,z_j\}
\end{align}
On the other hand, $|ij\rangle$ represents a chord stretching between the $i$-th and $j$-th strand of a given braid. We define a projection:
\begin{align}
\mathcal{A}^{(1)}(braid) &\xrightarrow{p_1} (\mathbb{C}^2-\Delta)/S_2 \\
|ij\rangle(Z) &\mapsto \{z_i,z_j\}
\end{align}
It follows that we must have $|ij\rangle \dlog(z_i-z_j) \in \mathcal{A}^{(1)}(\text{braid})\times_{X_2} \Omega^1 (\log \!\vartriangle \!\!\mathbb{C})$.\\

\section{The Kontsevich integral as a generator}
Any given link $L$ can be put in braid form \cite{A}, a geometric braid $\tilde{\gamma} \in \mathbb{C} \times I$ whose closure yields back the link we started with. Thus we regard links as being equivalent to their geometric braids. We regard a two strands geometric braid in $\mathbb{C} \times I$ with 4 boundary points, two of which are in the plane $\mathbb{C} \times \{0\}$, the other remaining two in the plane $\mathbb{C} \times \{1\}$, as the boundary of a ribbon of $\mathbb{C} \times I$ whose intersections with the planes $\mathbb{C} \times \{0\}$ and $\mathbb{C} \times \{1\}$ are exactly those 4 points. We refer to such a ribbon as being the ribbon associated to that particular geometric braid. Thus we regard two-strands geometric braids as being equivalent to their associated ribbons, since we can go from one to the other. Consequently, given a ribbon in $\mathbb{C} \times I$, we can also refer to its boundary strands as its two associated strands. We also regard a unique strand as being associated to a ribbon of vanishing width in an obvious way. This generalizes easily to geometric braids with $N$ strands; of the $(N-1)!$ associated ribbons connecting them, a smaller number is necessary to fully recover the $N$ strands of the braid, as well as to position the strands with respect to one another. In doing so we keep in mind that a geometric braid $\tilde{\gamma}$ is the lift of some loop $\gamma$ in configuration space $X_N$. To give the positioning of two strands by means of their associated ribbon is an equivalence relation, and by transitivity it follows that all we need is $(N-1)$ of those $(N-1)!$ ribbons. We now regard a given ribbon in $\mathbb{C} \times I$ associated to two strands of a geometric braid as the closure of the set of horizontal chords from one point of either of its associated strands to the point on its other associated strand. Such a closure defines a ruled surface whose underlying ribbon is none other than the ribbon we started with. This generalizes easily to the presence of $N-1$ ribbons. One may therefore think that to recover a link $L$ it is sufficient to have $N-1$ well-chosen ribbons, or equivalently the closures of $N-1$ sets of chords between their boundaries. What the Kontsevich integral does is a lot more. The Kontsevich integral is valued in the graded algebra $\mathcal{A}$ of chord diagrams. The metric aspect necessary for locating geometric objects with respect to one another is encoded in the local coordinates on $\mathbb{C} \times I$ for such objects. For instance two strands of a geometric braid are viewed as the lift of two paths in $X_2$ given by two functions $z_1(t)$ and $z_2(t)$, $t \in I$. Ultimately the Kontsevich integral is invariant under horizontal deformations, so what is of most interest to us is the winding of strands around one another, and thus we are led to considering not differences $z_1 - z_2$ but logarithms of such differences as those pick up crossings between strands. Thus a horizontal chord based at two points of two different strands along with the logarithm of the difference of the two complex variables locating these two points is sufficient. For $P$ an applicable pairing, $Z=\{z_1,z_2\}$ a point of $X_2$, $|P\rangle (Z)$ the chord between the two points $z_1$ and $z_2$, the object we are looking at is $(|P\rangle (Z), \frac{1}{2 \pi i}\log (z_1-z_2))$. Given a two strands geometric braid $\tilde{\gamma} = \{\tilde{\gamma}_1 , \tilde{\gamma}_2 \}$ with $\tilde{\gamma}_i$ the arc corresponding to the graph $\Gamma_i=\{(z_i (t), t) | t \in I \}$, obtained from lifting $\gamma_i = \{z_i(t) | t \in I \} $ in $X_2$, $i=1,2$, we regard its associated ribbon as the surface underlying the ruled surface obtained as the closure of the set $\{|P\rangle (z_1(t), z_2(t)) | t \in I \}$. If it is clear how to reproduce the ribbon, it is not clear however how to implement such a closure. A first step towards achieving this in Kontsevich integral computations is to fatten chords and to consider germs of chords based at small neighborhoods of points at which they are located, which are given as the intersections of small open balls in $\mathbb{C} \times I$ centered about those points with the strands on which the points are located. We denote by $\delta$ such an operation, and by $\delta Z$ such a neighborhood. We write:
\beq
\delta |P\rangle (Z)=|P\rangle (\delta Z)
\eeq
where $|P\rangle(Z)$ is a chord with support the point $Z$ while $|P\rangle(\delta Z)$ is the same chord with support in a neighborhood $\delta Z$ of the point $Z$, and with its feet located at $Z$. Once such an object is defined we can define differentials in such a neighborhood $\delta Z$, and thus we can consider the log differential $\dlog(z_1-z_2)$. This leads to considering densities defined as follows:
\beq
\delta (|P\rangle (Z), \log (z_1-z_2))= |P\rangle (\delta Z) \dlog (z_1-z_2)
\eeq
Further if $\overline{\cup _{Z \in \tilde{\gamma}}|P\rangle (Z)}$ does reproduce a ruled surface, we cannot easily incorporate the logarithms in such a closure. Summing over densities is possible however, and this is done via an integration. At the level of chord diagrams what used to be a simple union of chord diagrams based at a point can now be implemented by taking the concatenation of chord diagrams based in a neighborhood of points of $\tilde{\gamma}$, as such neighborhoods can be concatenated. This leads us to defining the following product. For two applicable pairings $P$ and $P'$ of degree one, $|P\rangle (Z_a)$ and $|P'\rangle (Z_b)$ based at two different points distant from one another, then we define:
\begin{align}
|P\rangle &(\delta Z_a) \dlog (z_{a,1} - z_{a,2}) \cdot |P'\rangle (\delta Z_b) \dlog (z_{b,1}-z_{b,2}) \nonumber \\
&=|P\rangle (\delta Z_a)|P'\rangle (\delta Z_b)\dlog (z_{a,1} - z_{a,2})\dlog (z_{b,1}-z_{b,2})
\end{align}
If the two points $Z_a$ and $Z_b$ are close together, then we can regard $\delta Z_a$ and $\delta Z_b$ as being essentially the same neighborhood and we define the product $|P\rangle (\delta Z_a) |P'\rangle(\delta Z_b)$ as being a concatenation strand-wise:
\beq
|P\rangle (\delta Z_a) |P'\rangle(\delta Z_b)=|P,P'\rangle (\delta Z_{\Lambda})
\eeq
$\Lambda$ being either of $a$ or $b$. In this situation, this leads to defining:
\begin{align}
|P\rangle &(\delta Z_a) \dlog (z_{a,1} - z_{a,2}) \cdot |P'\rangle (\delta Z_b) \dlog (z_{b,1}-z_{b,2})\nonumber \\
&=|P,P'\rangle (\delta Z_{\Lambda})\dlog (z_{a,1} - z_{a,2})\dlog (z_{b,1}-z_{b,2})
\end{align}
We generalize this easily to the product of more than two chord diagram valued log differentials. This gives rise to the non-commutative graded algebra $\delta \mathcal{A}(X_N \times I)$ with graded completion $\delta \overline{\mathcal{A}}(X_N \times I)$. Observe that the support of such chord diagrams was the braid itself in the original definition of the Kontsevich integral. Thus if we define $Z(t)=\tilde{\gamma} \cap \mathbb{C} \times \{t\}$, then this defines a fonction $Z$ on $I$. Then for $t$ fixed, $P$ fixed, $|P|=1$, the notation $|P\rangle (Z(t))$ makes sense. What we have is a minimal such definition where chords are not tangle chord diagrams per se \cite{ChDu} but rather are merely chord diagrams defined only locally. It is then easy to sum over such densities: for $m \geq 0$ fixed, for $P=(P_1,\cdots , P_m)$ fixed, we sum over all terms of the form $(\frac{1}{2 \pi i})^m\prod _{1 \leq i \leq m}|P_i\rangle (\delta Z(t_i)) \prod_{1 \leq i \leq m} \dlog (\vartriangle \!\! z [P_i] (Z(t_i))$ for $0 < t_1 < \cdots < t_m <0$. We then sum over all such choices of $P$'s for which $|P|=m$, and then finally sum over all $m \geq 0$:
\beq
\sum _{m \geq 0} \sum_{|P|=m} \int_{0 < t_1< \cdots < t_m <1} \frac{1}{(2 \pi i)^m}\prod _{1 \leq i \leq m}|P_i\rangle (\delta Z(t_i)) \prod_{1 \leq i \leq m} \dlog (\vartriangle \!\! z [P_i] (Z(t_i))
\eeq
We can easily generalize this to the case of a geometric braid with $N$ strands and we obtain the minimal Kontsevich integral $\Lambda$:
\beq
\Lambda=\sum _{m \geq 0} \sum_{|P|=m} \int_{0 < t_1< \cdots < t_m <1} \frac{1}{(2 \pi i)^m}\prod _{1 \leq i \leq m}|P_i\rangle (\delta Z(t_i)) \prod_{1 \leq i \leq m} \dlog (\vartriangle \!\! z [P_i] (Z(t_i))
\eeq
If we define:
\begin{align}
\Lambda_M&=\nonumber \\
&\sum _{0 \leq m \leq M} \sum_{|P|=m} \int_{0 < t_1< \cdots < t_m <1} \frac{1}{(2 \pi i)^m}\prod _{1 \leq i \leq m}|P_i\rangle (\delta Z(t_i)) \prod_{1 \leq i \leq m} \dlog (\vartriangle \!\! z [P_i] (Z(t_i))
\end{align}
then we can write:
\beq
\Lambda=\lim_{M \rightarrow \infty} \Lambda_M
\eeq
If we assume that the geometric braids we work with are smooth enough, then for $M$ large, $\Lambda_M$ is sufficient to geometrically produce ruled surfaces whose boundaries are the geometric braid we are seeking. We define the depth of a geometric braid to be the smallest value of $M$ for which we can recover the braid from studying the coefficients of $\Lambda_M$. From \cite{RG1} we know such a value is 1. The definition of depth will be most useful later when we generalize the Kontsevich integral to more complicated objects than simple geometric braids. Now however the minimal Kontsevich integral is not the integral as it was initially defined \cite{K}. Kontsevich used the skeleton as a support of the chord diagrams, and instead of considering chord diagrams defined locally, or equivalently germs of chord diagrams, one considers tangle chord diagrams with support on a geometric braid. This can easily be implemented by putting the chords of $\Lambda$ on the geometric braid $\tilde{\gamma}$. To do this we define an action of the graded completion of the non-commutative graded algebra $\oplus_{n \geq 0} \delta \mathcal{A}^{(n)}( X_N \times I) \times ( \Omega^1 (\log \vartriangle \!\mathbb{C}))^n$ on braids by recurrence. If $|P|=1$, $Z \in X_N \times \{t\}$, $t \in I$, $1 \leq i \neq j \leq N$ are given, $\tilde{\gamma}_h$, $h=i,j$ the strands of $\tilde{\gamma}$ on which $|P\rangle $ is supported, then we define $|P\rangle (\delta Z)  \cdot \tilde{\gamma}=|P\rangle (\tilde{\gamma}(t))$ if $\tilde{\gamma}(t)=Z$, $ \tilde{\gamma}$ otherwise. Thus here recovering the geometric braid is not the point of computing the Kontsevich integral. We have:
\begin{align}
\Lambda \cdot \tilde{\gamma}&= \nonumber \\
\sum _{m \geq 0} &\sum_{|P|=m} \int_{0 < t_1< \cdots < t_m <1} \frac{1}{(2 \pi i)^m}\prod _{1 \leq i \leq m}|P_i\rangle (\delta \tilde{\gamma}(t_i)) \prod_{1 \leq i \leq m} \dlog (\vartriangle \!\! z [P_i] (\tilde{\gamma}(t_i)) \cdot \tilde{\gamma}\\
=\sum _{m \geq 0} &\sum_{|P|=m} \int_{0 < t_1< \cdots < t_m <1} \frac{1}{(2 \pi i)^m}\prod _{1 \leq i \leq m}|P_i\rangle (\delta \tilde{\gamma}(t_i)) \cdot \tilde{\gamma}\prod_{1 \leq i \leq m} \dlog (\vartriangle \!\! z [P_i] (\tilde{\gamma}(t_i)) \\
=\sum _{m \geq 0} &\sum_{|P|=m} \int_{0 < t_1< \cdots < t_m <1} \frac{1}{(2 \pi i)^m} |P\rangle (\tilde{\gamma}T) \prod_{1 \leq i \leq m} \dlog (\vartriangle \!\! z [P_i] (\tilde{\gamma}(t_i))\\
&=Z(\tilde{\gamma})
\end{align}
where we have used the notation $T=(t_1,\cdots, t_m)$. This simplifies as follows:
\begin{align}
\sum _{m \geq 0} \sum_{|P|=m} &\int_{0 < t_1< \cdots < t_m <1} \frac{1}{(2 \pi i)^m} |P\rangle (\tilde{\gamma}T) \prod_{1 \leq i \leq m} \dlog (\vartriangle \!\! z [P_i] (\tilde{\gamma}(t_i))= \nonumber \\
&\sum _{m \geq 0} \sum_{|P|=m} |P\rangle (\tilde{\gamma}) \int_{0 < t_1< \cdots < t_m <1} \frac{1}{(2 \pi i)^m}  \prod_{1 \leq i \leq m} \dlog (\vartriangle \!\! z [P_i] (\tilde{\gamma}(t_i))
\end{align}
as all the representatives $|P\rangle (\tilde{\gamma}T)$ are elements of the homeomorphism class $|P\rangle (\tilde{\gamma})$ which we can factor out of the integral. This enables one to see that the Kontsevich integral is a sum over all degrees of chords, and for each degree a sum over all possible homeomorphism classes of chords of that particular degree, and for each such class a sum over all representatives, which is given by chords supported on the braid times the integral of an appropriate power of the log differentials which are none other than the densities necessary for performing such as sum. We can simplify this sum even further by defining an equivalence class on homeomorphism classes of tangle chord diagrams: define two pairings $P$ and $P'$ to be equivalent relative to $\tilde{\gamma}$ if one can go from one pairing to the other by sliding the chords of one pairing along $\tilde{\gamma}$ to obtain the chords of the other. If we close the geometric braid into a link, this is what we would obtain as the chords circle the link. Thus this equivalence relation becomes manifest once each tangle chord diagram is closed into a link. The resulting Kontsevich integral we denote by $\complement Z(\tilde{\gamma})$. We also write $\complement |P\rangle (\tilde{\gamma})=L_P$ if the geometric braid $\tilde{\gamma}$ closes into a link $L$. Then we can write:
\begin{align}
\complement Z(\tilde{\gamma})&=\complement \sum _{m \geq 0} \sum_{|P|=m} |P\rangle (\tilde{\gamma}) \int_{0 < t_1< \cdots < t_m <1} \frac{1}{(2 \pi i)^m}  \prod_{1 \leq i \leq m} \dlog (\vartriangle \!\! z [P_i] (\tilde{\gamma}(t_i))\\
&=\sum _{m \geq 0} \sum_{|P|=m} \complement |P\rangle (\tilde{\gamma}) \int_{0 < t_1< \cdots < t_m <1} \frac{1}{(2 \pi i)^m}  \prod_{1 \leq i \leq m} \dlog (\vartriangle \!\! z [P_i] (\tilde{\gamma}(t_i))\\
&=\sum _{m \geq 0} \sum_{|P|=m} L_P \int_{0 < t_1< \cdots < t_m <1} \frac{1}{(2 \pi i)^m}  \prod_{1 \leq i \leq m} \dlog (\vartriangle \!\! z [P_i] (\tilde{\gamma}(t_i))\\
&=\sum _{m \geq 0} \sum_{\substack{[P] \\ |P|=m}}\sum _{P' \in [P]} L_{P'} \int_{0 < t_1< \cdots < t_m <1} \frac{1}{(2 \pi i)^m}  \prod_{1 \leq i \leq m} \dlog (\vartriangle \!\! z [P'_i] (\tilde{\gamma}(t_i))\\
&=\sum _{m \geq 0} \sum_{\substack{[P] \\ |P|=m}} L_P \Big(\sum _{P' \in [P]}  \int_{0 < t_1< \cdots < t_m <1} \frac{1}{(2 \pi i)^m}  \prod_{1 \leq i \leq m} \dlog (\vartriangle \!\! z [P'_i] (\tilde{\gamma}(t_i)) \Big)
\end{align}

\section{The Kontsevich integral as a map from orbifolded cylinders}
We have the Kontsevich integral of geometric braids $\tilde{\gamma}$ embedded in $\mathbb{C} \times I$. One can of course consider the integral of links as well, the only addition being an overall sign for each chord diagram as follows:
\beq
Z(L)=\sum_{m \geq 0} \sum_{|P|=m} \frac{1}{(2 \pi i)^m} \int _{0 < t_1 < \cdots < t_m < 1} (-1)^{\epsilon (P)}L_P \prod_{ 1 \leq i \leq m} \dlog \vartriangle \!\! z[P_i](Z(t_i))
\eeq
and $Z(t)$ are local coordinates on $L$ and for each chord $\epsilon(P)$ counts the number of its feet ending on strands that are locally oriented down. We can easily generalize this integral to more general pictures such as oriented graphs whose vertices are univalent, trivalent (y or $\lambda$-shaped) or 4-valent (X-shaped). It is worth recalling at this point that the Kontsevich integral is defined for Morse links, and correspondingly we will assume that we do not have graphs with straight edges if we know that after rotation those edges may end up being horizontal. If $\{ G_i | 1 \leq i \leq q \}$ is a collection of graphs of $\mathbb{C} \times I$, we can compute $Z( \amalg _i G_i)$. In doing such a computation, we can study the behavior of $Z$ as graphs are moved in $\mathbb{C} \times I$, thereby introducing a time dependence $\tau$ in the computation of $Z$. Such a dynamic picture can be implemented as follows: for a graph $G$ and a path $\alpha$ in $\mathbb{C} \times I$ inducing a tangent vector field $X$, then moving $G$ along $\alpha$ means at time $\tau$ each point of $G$ moves in the direction given by the vector $X(\tau)$, and this for all $\tau \in [0,1]$. In other terms moving a graph along a path means moving it as a single block. Observe that $Z(G)$ of a single graph $G$ is time independent as the movement of the graph $G$ along any path will not alter $Z(G)$. As soon as we consider two or more graphs however, $Z(\amalg_{i \geq 2} G_i)$ becomes non-trivial as at least one of the graph is moved relative to the others. A trivial example is provided by two non-parallel strands with same orientation, with highest and lowest points in the same respective planes $\mathbb{C} \times \{2/3\}$ and $\mathbb{C} \times \{1/3\}$, with a separation of $a$ at the top, a separation of $b$ at the bottom. The degree one term of the Kontsevich integral of such a picture is $\log(a/b)/2 \pi i$, while if we move up either strand by a third of a unit, the resulting Kontsevich integral is trivial as there are no longer any chords between the two strands. We also consider the rotation of graphs with respect to a point. For a graph $G$ and a fixed point $p$ of $\mathbb{C} \times I$, we can rotate $G$ with respect to that point. The resulting Kontsevich integral will not be invariant under such a rotation as it is known that $Z$ depends on a choice of time axis \cite{K}. One can trivially convince oneself of this fact; the Kontsevich integral of the U-shaped unknot is non trivial, it is commonly denoted $\nu^{-1}$, whereas the Kontsevich integral of the same unknot rotated sideways by ninety degrees is trivial. Since we consider moving graphs, the point $p$ will not be fixed throughout but will change with time so we consider another path $\beta$ such that $\beta (\tau)$ will be the desired point at time $\tau$ with respect to which a graph is rotated.\\

Each graph $G_i$, $1 \leq i \leq q$ moves along a particular path $\alpha_i$ and has a particular curve of center of rotation points $\beta_i$. Thus what we have is a functional $Z(\amalg_i G_i)[\alpha_1][\beta_1][\cdots][\alpha_q][\beta_q]$ and we have that the Kontsevich integral appears as a map from $q$ cylinders $S^1 \times I$ to $\overline{\mathcal{A}}(\amalg_i G_i)$:
\beq
Z(\amalg_i G_i)[\alpha_1][\beta_1][\cdots][\alpha_q][\beta_q]:   \otimes ^q S^1 \times I \rightarrow \overline{\mathcal{A}}(\amalg_i G_i)
\eeq
Now as graphs move towards one another, the presence of logarithmic differentials in the expression for $Z$ may lead to singularities. As points connected by a chord get closer together the corresponding log differentials give rise to coefficients that are increasing in value. If contact occurs we distinguish two cases. If at the point of contact we have a vertex that is not y, $\lambda$ or X-shaped, we do have an infinite result. If the point of contact is trivalent (y or $\lambda$-shaped) or 4-valent (X-shaped), we have what we call a vanishing singularity for then the coefficient of such a resulting graph is finite by virtue of the framing independence relation. A first remark is that as the number of components is reduced we either have singularities or vanishing singularities, which points to the fact that the Kontsevich integral may be ultimately defined on a stratified space, something we will go into in a forthcoming paper. Observe that if we rotate graphs, what appears to be singularities may disappear altogether. For illustrative purposes, consider a circle and a strand at an angle that moves towards the circle and touches it say at the point of intersection of the horizontal line going through the center of the circle. This is a vanishing singularity. If however we move this strand around the circle in such a manner that it touches the circle on the vertical line going through the center of the circle, then there never was a singularity. \\

Thus two graphs $G_1$ and $G_2$ can be brought into contact:
\begin{itemize}
\item[-] At some points $(\sigma_1, \tau_1)$ and $(\sigma_2, \tau_2)$ and are therefore identified. Thus $Z$ is defined on $S^1 \times I \otimes S^1 \times I/\{(\sigma_1, \tau_1) \sim (\sigma_2, \tau_2)\}$ where $Z(G_1 \amalg G_2)$ is defined away from the singular point on the quotient and $Z(G_1 \cup G_2)$ is defined exactly at the point where the two graphs are brought into contact to form what we call $G_1 \cup G_2$. If this results in a singularity of $Z$, we mark this point by an ``X", a point otherwise. The point of contact depends on the choice of $G=\{G_1, \cdots, G_q\}$, $\alpha =\{ \alpha_1, \cdots, \alpha_q \}$ and $\beta=\{ \beta_1, \cdots, \beta_q \}$, thus we will denote by $\sim_{G \alpha \beta}$ such an identification and by $S^1 \times I \otimes S^1 \times I /\sim_{G \alpha \beta}$ the resulting space on which $Z$ is defined. In so doing we adopt the Knot Theory point of view that tensor products can be represented by objects side by side. In this manner the identification can be easily visualized as being a simple gluing between cylinders, leading to a singular space that we will refer to as an identifold. We will refer to those glued cylinders as id-folded cylinders for short.
    \item[-] The two graphs $G_1$ and $G_2$ can be brought into contact and $G_1 \cup G_2$ exists along some path $(\alpha, \beta)$ which can either be given by $(\alpha_1, \beta_1)$ or $(\alpha_2, \beta_2)$. This corresponds to some values $(\sigma_i, \tau_i) \in [\theta_{i1}, \theta_{i2}]\times [a_i, b_i]$, $i=1,2$ on their respective cylinders being identified, leading to a common arc. Such an identification is taken into account by saying that $Z$ is a map on $S^1 \times I \otimes S^1 \times I /\sim_{G \alpha \beta}$. Subarcs of the arc of contact are drawn as a solid line if along such subarcs $Z$ is singular. Subarcs on which $Z$ is well-defined are drawn as a dashed line. Now along the arc of contact, $G_1 \cup G_2$ may be brought into contact with other graphs. The two above steps can then be repeated, leading to a second identification of points or subarcs of this arc with points from a third cylinder. All of this is still taken into account by working with the quotient $\otimes ^3 S^1 \times I/ \sim_{G \alpha \beta}$.
    \item[-] The two graphs $G_1$ and $G_2$ are brought into contact along some area. One instance where this happens is in the event that we have two circles of radius $0.5$ units centered at $(0,0)$ and $(1,0)$ respectively, each moving straight up, the circle on the left rotating counterclockwise as it moves, the one on the right rotating clockwise. Those two graphs are in contact for all times and angles. In that situation the two cylinders corresponding to those two circles are identified.
    In the event that contact occurs only for areas $\Sigma_1$ and $\Sigma_2$ possibly ending on either or both boundaries of $S^1 \times I$, we identify such areas to yield a common area $\Sigma$ in $S^1 \times I \otimes S^1 \times I / \sim_{G \alpha \beta}$. Subareas of $\Sigma$ over which contact between $G_1$ and $G_2$ results in a singularity for $Z$ are delimited by a solid line, a dashed line otherwise.
    The resulting graph $G_1 \cup G_2$ can further be brought into contact with other graphs, resulting in the area of contact in $\otimes^3 S^1 \times I / \sim_{G \alpha \beta}$ having points, arcs or subareas being identified with points from a third cylinder.
\end{itemize}
This has been done for two or three graphs being brought into contact but can easily be generalized to $q$ graphs $G_i$ being brought into contact, the geometry of contact still being taken into account by working with $\otimes ^q S^1 \times I / \sim_{G \alpha \beta}$, on which $Z$ is defined. We denote such an identifold by $IdX(S^1 \times I, G, \alpha, \beta)$ and by $IdX_{S^1 \times I}$ the set of all such identifolds.

More generally, if $\Gamma(\mathbb{C} \times I)$ denotes the set of graphs embedded in $\mathbb{C} \times I$, then $Z$ is an element of $F(  \Gamma (\mathbb{C} \times I)^q, F( (P ( \mathbb{C} \times I ))^{2q}, F( IdX_{S^1 \times I}, \overline{\mathcal{A}}( \Gamma ( \mathbb{C} \times I))^q)))$. In stages:
\begin{align}
Z: \Gamma (\mathbb{C} \times I)^q &\rightarrow F( (P ( \mathbb{C} \times I ))^{2q}, F(IdX_{S^1 \times I}, \overline{\mathcal{A}}( \Gamma ( \mathbb{C} \times I))^q))  \nonumber \\
\amalg_{1 \leq i \leq q} G_i &\mapsto Z(\amalg_{1 \leq i \leq q} G_i)
\end{align}
To those $q$ graphs in $\mathbb{C} \times I$, we can associate curves for translations as well as curves for rotations as follows:
\begin{align}
Z(\amalg_{1 \leq i \leq q} G_i): P ( \mathbb{C} \times I ))^{2q} \rightarrow F(IdX_{S^1 \times I}, \overline{\mathcal{A}}( \Gamma ( \mathbb{C} \times I))^q) \nonumber \\
\times_{1 \leq i \leq q} (\alpha_i, \beta_i) \mapsto Z(\amalg_i G_i, \times_i (\alpha_i, \beta_i))
\end{align}
Once these paths are defined, the resulting Kontsevich integral can be seen as being a map of towers to the graded algebra of chord diagrams with support the $q$ diagrams that were initially chosen:
\begin{align}
Z(\amalg_i G_i, \times_i (\alpha_i, \beta_i)): IdX_{S^1 \times I} \rightarrow \overline{\mathcal{A}}(\amalg_q G_i) \nonumber \\
\otimes^q S^1 \times I / \sim_{G, \alpha \beta} \mapsto Z(\amalg_i G_i, \times^q (\alpha_i, \beta_i))
\end{align}
At the second stage above, we can study the deformations of $Z$ under deformations in the space of paths in $\mathbb{C}\times I$. For a graph $G$ in $\Gamma ( \mathbb{C} \times I)$ which is the result of glueing $N(G)$ arcs $\tilde{\gamma}_i$ in $\mathbb{C} \times I$, with associated tangent vectors $X_i$, $ 1 \leq i \leq N(G)$, then we define the tangent space to $\Gamma(\mathbb{C} \times I)$ at $G$ to be given by:
\beq
T_G \Gamma(\mathbb{C} \times I)=\{X_i | 1 \leq i \leq N(G) \}
\eeq
Such arcs $\tilde{\gamma}_i$ are lifts of paths $\gamma_i$ in the complex plane and thus are given by $(\gamma_i(t), t)$ for $t \in I$, which we denote by $\tilde{\gamma_i}(t)$. We denote by $d \tilde{\gamma}_i$ a differential along such arcs $\tilde{\gamma}_i$, duals to the vectors $X_i$. Then the cotangent space to a graph $G$ is defined to be:
\beq
T^*_G \Gamma(\mathbb{C} \times I) = \{d \tilde{\gamma}_i | 1 \leq i \leq N(G) \}/{G' \amalg G''=G}
\eeq
We write:
\beq
\delta G =\sum_{1 \leq i \leq N(G)}\lambda_i d\tilde{\gamma}_i
\eeq
the formal deformation of $G$ where the $\lambda_i$'s are coefficients. We can also deform the paths $\alpha$ and $\beta$ which leads to defining the tangent space:
\begin{align}
T_{\alpha, \beta}P(\mathbb{C}\times I)^q &=\{\delta(\alpha, \beta)\} \\
&=\{((d\alpha_1, d\beta_1), \cdots, (d\alpha_q, d\beta_q))\}
\end{align}
All such deformations induce deformations of $\otimes^q S^1 \times I / \sim_{G \alpha \beta}$. Observe that absent any knowledge of $G$, $\alpha$ or $\beta$, knowing this quotient space we can determine when divergences for $Z$ arise, thereby presenting such identifolds as a blueprint for studying the singularities of the Kontsevich integral.

\end{document}